\documentclass[12pt]{article}
\usepackage{amsmath}
\usepackage{amssymb}
\newcommand{\poir}{P[f](re^{i\theta})}
\newcommand{\hk}{{\cal HK}}
\newcommand{\bv}{{\cal BV}}
\newcommand{\poi}{P[f]}

\newcommand{\R}{{\mathbb R}}
\newcommand{\fn}{\!:\!}
\newcommand{\lsum}{\sum\limits}
\newcommand{\llim}{\lim\limits}
\newcommand{\lint}{\int\limits}
\newcommand{\lsup}{\sup\limits}
\newcommand{\linf}{\inf\limits}
\newcommand{\qed}{\mbox{$\quad\blacksquare$}}
\newcommand{\thn}{\theta_n}
\newcommand{\ret}{re^{i\theta}}
\newtheorem{theorem}{Theorem}

\newtheorem{corollary}[theorem]{Corollary}
\newtheorem{remarks}[theorem]{Remarks}
\newtheorem{example}[theorem]{Example}
\begin{document}
\hspace{-2cm}
\raisebox{12ex}[1ex]{\fbox{{\footnotesize
To appear in {\it Canadian  Mathematical Bulletin.}
Accepted December 24, 2003.
}}}
$~$ \\ [.5in]

\begin{center}
{\large\bf Estimates of Henstock--Kurzweil Poisson integrals}
\vskip.25in
Erik Talvila\footnote{Research partially supported by the
Natural Sciences and Engineering Research Council of Canada.
An adjunct appointment in the Department of Mathematical and
Statistical Sciences, University of Alberta, made valuable library
and computer resources available.}\\ [2mm]
{\footnotesize
Department of Mathematics and Statistics \\
University College of the Fraser Valley\\
Abbotsford, BC Canada V2S 7M8\\
Erik.Talvila@ucfv.ca}
\end{center}

{\footnotesize
\noindent
{\bf Abstract.} If $f$ is a real-valued function on $[-\pi,\pi]$ that
is Henstock--Kurzweil integrable, let $u_r(\theta)$ be its Poisson
integral.  It is shown that $\|u_r\|_p=o(1/(1-r))$ as $r\to 1$
and this estimate is sharp for $1\leq p\leq\infty$.  
If $\mu$ is a finite Borel measure and $u_r(\theta)$ is its Poisson
integral then for each $1\leq p\leq \infty$ the estimate
$\|u_r\|_p=O((1-r)^{1/p-1})$ as $r\to 1$ is sharp.
The Alexiewicz
norm estimates $\|u_r\|\leq\|f\|$ ($0\leq r<1$) and $\|u_r-f\|\to 0$
($r\to 1$) hold. These estimates lead to two uniqueness theorems for
the Dirichlet problem
in the unit disc with Henstock--Kurzweil integrable boundary data.
There are similar growth estimates when $u$ is in the harmonic Hardy
space associated with the Alexiewicz
norm and when $f$ is of bounded variation.
\\
{\bf 2000 subject classification:} 26A39, 31A20
}
\section{Introduction}
In this paper we consider estimates of Poisson integrals on the
unit circle
with respect to Alexiewicz and $L^p$ norms.
Define the open disk in $\R^2$ as
$D:=\{re^{i\theta}\,|\,0\leq r<1,\, -\pi\leq\theta< \pi\}$.
The Poisson kernel is $\Phi_r(\theta):=(1-r^2)/[2\pi(1-2r\cos\theta+r^2)]
=[1+2\sum_{n=1}^\infty r^n\cos(n\theta)]/(2\pi)$.  
Let $f\fn \R\to\R$ be $2\pi$-periodic.
The
Poisson integral of $f$ is its convolution with the Poisson kernel
$$
\poir:=f\ast\Phi_r(\theta)=\int_{-\pi}^\pi f(\phi)\Phi_r(\phi-\theta)\,d\phi.
\label{poissonintegral}
$$
Since $\partial D$ has no end points, 
an appropriate form of  
the Alexiewicz norm of $f$ is $\|f\|:=\sup_{I\subset \R}\left|
\int_{I}f\right|$ where $I$ is an interval in $\R$ of length not exceeding
$2\pi$.  
Let $\hk$ denote the $2\pi$-periodic 
functions 
$f\fn\R\to\R$ with finite Alexiewicz norm. 
Of course, with the same periodicity convention, 
$L^p\subsetneq\hk$ for all $1\leq p\leq\infty$.
Write $\|f\|_A$ for the Alexiewicz norm
over set $A$.
The Alexiewicz norm
is discussed in \cite{swartz2}.
The variation of $f$ over one period  is denoted 
$Vf$.  The set of $2\pi$-periodic functions with finite
variation over one period is denoted $\bv$.
For a function $u\fn D\to \R$ we write $u_r(\theta)=u(\ret)$.

The Dirichlet problem, of finding a function harmonic in the disc 
with prescribed boundary values, is one of the foundational 
problems in elliptic partial differential equations.
An understanding of its solution has been a stepping stone
to the study of analytic functions in the complex plane 
and of the solutions of more general elliptic equations.
Due to the simple geometry of the disc there is an explicit
integral representation for solutions through \eqref{poissonintegral}.
As a Lebesgue integral, the Poisson integral has been studied intensively.
For the major results,   see
\cite{axler} and \cite{wolf}.

The following results are well known \cite{axler}.  Suppose
that $1\leq p \leq \infty$ and $f\in L^p$.
If $|\theta_0|\leq\pi$ and $z\in D$,
we say that $z\to e^{i\theta_0}$  nontangentially
if there is $0\leq \alpha
<\pi/2$ such that $z\to e^{i\theta_0}$ within the sector
$\{\zeta\in D:|\arg(\zeta-e^{i\theta_0})
|<\alpha\}$.
Write $u_r(\theta)=\poir$.  Then
\begin{align}
&u_r \text{ is harmonic in } D \label{1}\\
&\|u_r\|_p\leq\|f\|_p \text{ for all } 0\leq r <1.\label{2}\\
&\text{If } 1\leq p < \infty \text{ then } \|u_r-f\|_p\to 0  \text{ as }
r\to 1\label{3}\\
&u(\ret)\to f(\theta_0)  \text{ for almost all } \theta_0 \text{ as } z\to  
e^{i\theta_0} \text{ nontangentially
in } D.\label{4}
\end{align}

We examine analogues of these results when
$f$ is Henstock--Kurzweil integrable (Theorem~\ref{theorem2}).  
We also prove that the growth
estimate $\|u_r\|_p=o(1/(1-r))$ is sharp for $f\in\hk$  and
$1\leq p\leq \infty$ (Theorem~\ref{theorem1}).
If $\mu$ is a finite Borel measure and $u_r(\theta)$ is its Poisson
integral then for each $1\leq p\leq \infty$ the estimate
$\|u_r\|_p=O((1-r)^{1/p-1})$ as $r\to 1$ is sharp (Remarks~\ref{remarks1}).
The Poisson integral 
of a function in $\hk$ need not be the difference
of two positive harmonic functions (Remarks~\ref{rem1.5}).
There are similar
growth estimates when $u$ is in $h^\hk$,  the
harmonic Hardy space associated with the Alexiewicz norm 
(Theorem~\ref{hardy}).  The Poisson integral provides an isometry
from $\hk$ into (but not onto) $h^\hk$ (Theorem~\ref{isometry}).
In Theorem~\ref{bv} we consider the above results for functions of
bounded variation.
Theorem~\ref{theorem4} and Theorem~\ref{theorem5} establish
uniqueness conditions for the Dirichlet problem using the Alexiewicz
norm.   Example~\ref{example} shows the applicability of the uniqueness
theorems.
All the results also hold when we use
the wide Denjoy integral \cite{celidze}.

Since $\Phi_r$ and $1/\Phi_r$ are of bounded variation on $\partial D$,
necessary and sufficient for the
existence of $\poi$ in $D$ is that $f$ be integrable, 
i.e., the Henstock--Kurzweil integral
$\int_{-\pi}^{\pi}\!f$ is finite.
	In \cite{benedicks},
integration by parts was used to show that we can differentiate under
the integral sign.  This in turn shows that
$\poi$ is harmonic in $D$ and that $\poi\to f$ nontangentially,
almost everywhere on $\partial D$.
In \cite{celidze} (Theorem~4, p.~238),
necessary 
and sufficient conditions were given for
determining when a function that is harmonic in $D$ is the Poisson integral of
an $\hk$ function. Corresponding results when $\|u_r\|_p$
are uniformly bounded have been known for some time (for example, \cite{axler},
Theorem~6.13).

\section{Growth estimates} 
Our first result is to
show that for $1\leq p\leq \infty$, we have
$\|u_r\|_p=o(1/(1-r))$ and this estimate is sharp.
That is, 
$(1-r)\|u_r\|_p\to 0$ as $r\to 1$ ($1\leq p<\infty$) and 
$\sup_{\theta\in [-\pi,\pi]}(1-r)|\poir|\to 0$ as $r\to 1$ ($p=\infty$).
Thus, for $p=\infty$, the manner of approach to the boundary is unrestricted.
This same estimate for $p=\infty$ was obtained for
$L^1$ functions in \cite{wolf}.
We show these estimates are the best possible under our minimal
existence hypothesis.  The proof uses the inequality
\begin{equation}
\left|\int_{-\pi}^\pi fg\right|\leq 
\|f\|\left(\linf_{[-\pi,\pi]}|g|+Vg\right),\label{4.5}
\end{equation}
which is valid for all $f\in\hk$ and $g$ of bounded variation
on $[-\pi,\pi]$.  This was proved in \cite[Lemma~24]{talvilafourier}.

\begin{theorem}\label{theorem1}
Let $f\in\hk$.
For $\ret\in D$ let
$u_r(\theta)=P[f](\ret)$.
\begin{enumerate}
\item[{\rm (a)}] We have 
$\sup_{\theta\in[-\pi,\pi]}|\poir|=o(1/(1-r))$ as $r\to 1$
and this estimate is sharp
in the sense that if $\psi\fn D\to(0,\infty)$ and 
$\psi(\ret)=o(1/(1-r))$ as $r\to 1$ then there is a function $f\in\hk$
such that $P[f]\not=o(\psi)$ as $r\to 1$.
\item[{\rm (b)}] Let $1\leq p<\infty$.  
Then $\|u_r\|_p=o(1/(1-r))$ as $r\to 1$ and this estimate is sharp
in the sense that if $\psi\fn [0,1)\to(0,\infty)$ and
$\psi(r)=o(1/(1-r))$ as $r\to1$ then there is a function $f\in\hk$
such that $\|u_r\|_p\not=o(\psi(r))$ as $r\to1$.
\end{enumerate}
\end{theorem}

\noindent
{\bf Proof}:  (a)
Let $\Psi_r(\phi):=(1-r)^2/(1-2r\cos\phi+r^2)$ with 
$\Psi_1(0):=1$.  Let $0<\delta<\pi$.
Then
$$
\frac{2\pi(1-r)\poir}{1+r}=\lint_{|\phi-\theta|<\delta}f(\phi)\Psi_r(\phi-\theta)\,d\phi
\quad+\!\!\!\lint_{\delta<|\phi-\theta|<\pi}\!\!\!\!f(\phi)\Psi_r(\phi-\theta)\,d\phi.
$$
Given $\epsilon>0$, take $\delta$ small enough so that
$\|f\|_{[\theta-\delta,\theta+\delta]}<\epsilon$ for all $\theta$.
Using \eqref{4.5},
$$
\left|\,\lint_{|\phi-\theta|<\delta}f(\phi)\Psi_r(\phi-\theta)\,d\phi\right| 
\leq  
2\|f\|_{[\theta-\delta,\theta+\delta]}.
$$
And,
\begin{eqnarray*}
\left|\,\lint_{\theta+\delta}^{\theta-\delta+2\pi} f(\phi)\Psi_r(\phi-\theta)\,d\phi\right| & \leq & 
\|f\|\left[\frac{2(1-r)^2}{1-2r\cos\delta+r^2}-
\frac{(1-r)^2}{(1+r)^2}\right]\\
  & \to & 0 \quad\text{as } r\to 1.
  \end{eqnarray*}

To prove this estimate is sharp, suppose $\psi\fn D\to (0,\infty)$ is given.  
It suffices
to show that $\poi(r_ne^{i\theta_n})\not=o(\psi(r_ne^{i\theta_n}))$ for some
sequence $\{r_ne^{i\theta_n}\}\in D$ with $r_n\to 1$.
Take 
$0<\theta_n<\pi/2$ and decreasing to $0$.
Let $a_n=\psi(r_ne^{i\theta_n})$.
Take $0<\alpha_n\leq
\min(\pi/2, (\theta_{n-1}-\thn)/2,(\thn-\theta_{n+1})/2, 1-r_n)$ with
$\theta_0:=\pi$.
Then the intervals $(\thn-\alpha_n,\thn+\alpha_n)$ are disjoint and
$\cos(\alpha_n)\geq 1-\alpha_n^2/2$.
Let $f_n=\pi (1-r_n)a_n/\alpha_n$.  Define
$$
f(\phi)=\left\{\begin{array}{cl}
f_n, &|\phi-\thn|< \alpha_n\quad \text{for some } n\\
0, & \text{ otherwise.}
\end{array}
\right.
$$
Now,
\begin{eqnarray*}
2\pi\poi(r_ne^{i\theta_n}) & = & (1-r_n^2)\lsum_{k=1}^\infty
f_k\int_{\theta_k-\alpha_k}^{\theta_k+\alpha_k}
\frac{d\phi}{
r_n^2-2r_n\cos(\thn-\phi)+1}\\
 & \geq & \frac{2(1-r_n^2)f_n\,\alpha_n}{r_n^2-2r_n\cos(\alpha_n)+1}\\
 & \geq & \frac{2(1+r_n)(1-r_n)f_n\,\alpha_n}{(1-r_n)^2+r_n\alpha_n^2}.
\end{eqnarray*}
Hence,
$\poi(r_ne^{i\theta_n})
\geq a_n$.  And, 
$f\in L^1$ if $\sum f_k\,\alpha_k =\pi\sum (1-r_k)a_k <\infty$.
Since $(1-r_k)a_k\to 0$ there is a subsequence $\{(1-r_n)a_n\}_{n\in I}$ 
defined by
an unbounded index set $I\subset {\mathbb N}$ such that $\sum_{k\in I}(1-r_k)
 a_k 
<\infty$.  
Now take $f(\phi)=f_n$ when $|\phi-\theta_n|<\alpha_n$ for some $n\in I$ and
$f(\phi)=0$, otherwise. 
Then, $f\in L^1$ and 
$\poi(r_ne^{i\theta_n})\geq \psi(r_ne^{i\theta_n})$
for all $n\in I$.

(b) 
Suppose $1\leq p< \infty$.
From part (a), we can write $u_r(\theta)=w_r(\theta)/(1-r)$
where $\sup_{\theta\in [-\pi,\pi]}|w_r(\theta)|\to 0$ as $r\to 1$.
And, $w_r$ is periodic and
real analytic on $\R$ for each $0\leq r<1$.  Let $1\leq p<\infty$.  Then
\begin{eqnarray*}
\|u_r\|_p & = & \frac{1}{1-r}\left[\,\,\int_{-\pi}^{\pi}|w_r(\theta)|^p\,d\theta
\right]^{1/p}\\
 & \leq & \frac{(2\pi)^{1/p}}{1-r}\lsup_{\theta\in [-\pi,\pi]}|w_r(\theta)|.
 \end{eqnarray*}
 Hence, $\|u_r\|_p=o(1/(1-r))$ as $r\to 1$.

To prove this estimate is sharp, first consider $p=1$.  Let $\psi\fn
[0,1)\to(0,\infty)$ with $\psi(r)=o(1/(1-r))$ be given.
Although $\hk$ is not complete it is barrelled \cite{swartz2}.
The Uniform Boundedness Principle \cite{swartz1} applies and this
shows the existence of $f\in\hk$ such that $\|u_r\|_1\not=o(\psi(r))$.
We can see this as follows.

Define $r_n=1-1/n$ for $n\in{\mathbb N}$.  Let 
$f_n(\theta)=\psi(r_n)\sin(n\theta)$.  Then
\begin{eqnarray*}
\|f_n\| & = & \psi(r_n)\int_{0}^{\pi/n}\sin(n\theta)\,d\theta\\
 & = & 2\,\psi(r_n)/n\\
 & = & 2(1-r_n)\psi(r_n)\to 0\text{ as } n\to\infty.
\end{eqnarray*}

For $0\leq r<1$, define $S_r\fn\hk\to L^1$ by $S_r[f](\theta)=
P[f](\ret)/\psi(r)$ for each $f\in\hk$.  Write
$u_r(\theta)=P[f](\ret)$.  Using \eqref{4.5},
\begin{eqnarray}
\|u_r\|_1 & = & \int_{-\pi}^\pi \left|
\int_{-\pi}^\pi f(\phi) \Phi_r(\phi-\theta)\,d\phi\right|d\theta\notag\\
 & \leq & 2\pi\|f\|\left[\inf\Phi_r+V\Phi_r\right]\notag\\
 & = & \|f\|\left(\frac{1+6r+r^2}{1-r^2}\right).\label{6}
\end{eqnarray}
Therefore, $\|S_r\|\leq \frac{1+6r+r^2}{\psi(r)(1-r^2)}$ and, for each 
$0\leq r<1$, $S_r$ is a bounded linear operator from $\hk$ to $L^1$.

We have $S_r[f_n](\theta)=\psi(r_n)r^n\sin(n\theta)/\psi(r)$ so that
\begin{eqnarray}
\|S_{r_n}[f_n]\|_1 & = & r_n^n\int_{-\pi}^\pi |\sin(n\theta)|d\theta\notag\\
 & = & 4\left(1-1/n\right)^n\label{7}\\
 & \to & 4/e\quad \text{as } n\to\infty.\notag
\end{eqnarray}
It follows that $\{S_{r_n}\}$ is not equicontinuous.  The Uniform 
Boundedness Principle \cite[Theorem 11, p.~299]{swartz1} now shows
that $\{S_{r_n}\}$ is not pointwise bounded on $\hk$.  Hence, there
exists $f\in\hk$ such that $\sup_{n}\|u_{r_n}\|_1/\psi(r_n)=\infty$
and hence $\|u_r\|_1\not=o(\psi(r))$ as $r\to 1$.

The case $p>1$ is similar.  In place of \eqref{6}, we have
$\|u_r\|_p\leq (2\pi)^{1/p-1}\|f\|(1+6r+r^2)/(1-r^2)$.
And, in place of \eqref{7}, 
$$
\|S_{r_n}[f_n]\|_p=\left(1-1/n\right)^n
\left[\frac{2\sqrt{\pi}\,\Gamma((1+p)/2)}{\Gamma(1+p/2)}\right]^{1/p}.
\qed
$$

\begin{remarks}\label{remarks1}
{\rm 
The little {\it oh} order relation in Theorem~\ref{theorem1}(a) is
false for measures.  If $\mu$ is a finite
Borel measure on $[-\pi,\pi)$, write $u_r(\theta)=
P[\mu](\ret)$. Then $\|u_r\|_\infty\leq \Phi_r(0)
\mu([-\pi,\pi))=O(1/(1-r))$. 
The Dirac measure shows 
this estimate
is sharp.

For $1\leq p<\infty$, let $u_r(\theta)=P[\mu](\ret)$.  The
Minkowski inequality for integrals \cite[Theorem~6.19]{folland} gives
\begin{eqnarray*}
\|u_r\|_p & = & \left\|\int_{-\pi}^\pi \Phi_r(\phi-\cdot)\,d\mu(\phi)\right
\|_p\\
 & \leq & \int_{-\pi}^\pi \|\Phi_r(\phi-\cdot)\|_p\,d\mu(\phi)\\
 & = & \|\Phi_r\|_p\,\mu([-\pi,\pi)).
\end{eqnarray*}
And, for $\mu=\delta$, the Dirac measure, let $v_r(\theta)=P[\delta](\ret)$.
Then
\begin{eqnarray}
\|v_r\|_p & = & \|\Phi_r\|_p\notag\\
 & = & \frac{1-r^2}{2\pi}\left(\int_{-\pi}^\pi \frac{d\phi}{
(1-2r\cos\phi+r^2)^p}\right)^{1/p}\notag\\
& = & (2\pi)^{1/p-1}(1-r^2)^{1/p-1}\,\left[_2F_1(1-p,1-p;1;r^2)\right]^{1/p}
.\label{hyper}
\end{eqnarray}
Line \eqref{hyper} is from integral 3.665.2 in \cite{gradshteyn}
and the hypergeometric linear transformation \cite[9.131.1]{gradshteyn}.
For these values of the parameters, the hypergeometric function is
bounded for $0\leq r\leq 1$ and $_2F_1(1-p,1-p;1;1)=\Gamma(2p-1)/\Gamma^2(p)
\not=0$.  It follows that $\|u_r\|_p=O((1-r)^{1/p-1})$ as $r\to 1$.  The
Dirac measure shows this estimate is sharp.

The estimate for $p=1$ appears as Theorem~6.4(a) in \cite{axler}.
\qed
}
\end{remarks}

Several results follow immediately from these estimates.
For $1\leq p<\infty$, denote the harmonic Hardy spaces by 
$h^p:=\{u\fn D\to\R\mid
\Delta u=0 \text{ in } D, \sup_{0\leq r<1}\|u_r\|_p<\infty\}$.
And, $h^\infty$ is the set of bounded functions that are
harmonic in $D$.
The harmonic Hardy space associated with the Alexiewicz norm is defined
$$
h^{\hk}:=\{u\fn D\to\R\mid\Delta u=0, \lsup_{0\leq r<1}\|u_r\|<\infty\}.
$$
This is a normed linear space under the norm $\|u\|_{\hk}:= \sup_{0\leq r<1}
\|u_r\|$.

\begin{corollary}\label{cor4}
For $1\leq p\leq\infty$ we have $h^p\subsetneq h^\hk$.
\end{corollary}

\noindent
{\bf Proof:} We have $h^q\subset h^p\subset h^{\hk}$ 
for all $1\leq p < q \leq \infty$.
And, by Theorem~\ref{theorem1}(b), there is $f\in\hk$ with $u_r(\theta):=P[f](\ret)$ and $\|u_r\|_1
\not=O(1)$. \qed 

\begin{remarks}\label{rem1.5}
{\rm 
There  is a function $f\in\hk$ such that $P[f]$ is not the difference of
two positive harmonic functions.
This follows since
functions in $h^1$ are characterised as being the difference of
two positive harmonic functions.  See \cite[Exercise~6.9]{axler}. \qed
}
\end{remarks}

When 
$u\in h^\hk$ we can get slightly  different estimates than in
Theorem~\ref{theorem1}.
(cf. 
\cite[Proposition~6.16 and Exercise~6.11]{axler}).
\begin{theorem}\label{hardy}
Let $1\leq p\leq \infty$.
If $u\in h^{\hk}$ then $\|u_r\|_p\leq (2\pi)^{1/p}\frac{2r \|u\|_{\hk}}
{\pi(1-r)}$ for $1/2\leq r<1$ and 
$\|u_r\|_p\leq (2\pi)^{1/p}\frac{2 \|u\|_{\hk}}
{\pi}$ for $0\leq r\leq 1/2$.  (Replace the term $(2\pi)^{1/p}$ by
$1$ when $p=\infty$.)  The order relations are sharp as $r\to 1$.
\end{theorem}

\noindent
{\bf Proof:} Fix $z=\ret\in D$ and $0<t<1-r$.  If $0<t\leq r$ then,
by the Mean Value Property for
harmonic functions,
$u(z) 
  =  (\pi t^2)^{-1}\int_{r-t}^{r+t}\int_{\theta-\theta_0}^{\theta+
\theta_0}u(\rho e^{i\phi})\,d\phi\,\rho\, d\rho$,
where $\theta_0=\arccos[(r^2+\rho^2-t^2)/(2r\rho)]$ and $0\leq\theta_0
\leq\pi/2$.
Hence,
$$
|u(z)|  \leq  
 \frac{1}{\pi t^2}\int_{r-t}^{r+t}\rho\,d\rho\lsup_{|\rho-r|<t}
\left|\int_{\theta-\theta_0}^{\theta+
\theta_0}u(\rho e^{i\phi})\,d\phi\right|
  \leq  \frac{2r}{\pi t}\|u\|_{\hk}.
$$
Now let $t \to 1-r$ when $1/2\leq r<1$ and let $t\to r$ when 
$0\leq r\leq 1/2$.  This establishes the estimates for $p=\infty$.
The estimates for $1\leq p<\infty$ follow from this.
The case $r=0$ is similar.

Note that if $u(\ret)=\Phi_r(\theta)$ then $\|u\|_{\hk}=1$
and $\|\Phi_r\|_\infty=(1+r)/[2\pi(1-r)]$.
So, the order relation for $\|u_r\|_\infty$ is sharp as $r\to 1$.
For $1\leq p<\infty$, the implied order relation $O(1/(1-r))$ is
sharp as $r\to 1$ due to the example in the proof of 
Theorem~\ref{theorem1}(b).  For, suppose we are
given $\psi\fn[0,1)\to(0,\infty)$
with $\psi(r)=o((1-r)^{-1})$ as $r\to 1$.  From Theorem~\ref{theorem1}(b)
we know there is a function $f\in\hk$ such that if $u_r(\theta)=P[f](\ret)$
then $\limsup_{r\to 1}\|u_r\|_p/\psi(r)=\infty$.  And,
by the following Theorem~\ref{theorem2}(a), $\|u_r\|\leq\|f\|$ so 
$u\in h^{\hk}$.
\qed

Now consider the analogues of 
\eqref{2} and \eqref{3}
for the 
Alexiewicz norm.

\begin{theorem}\label{theorem2}
Let $f\in\hk$.  For $\ret\in D$ define
$u_r(\theta):=P[f](\ret)$. Then
\begin{enumerate}
\item[{\rm (a)}] $\|u_r\|\leq \|f\|$ for all $0\leq r <1$, i.e., 
$\|u\|_{\hk}\leq \|f\|$.
\item[{\rm (b)}] $\|u_r-f\|\to 0$ as $r\to1$
\item[{\rm (c)}] In (b), the decay of $\|u_r-f\|$  to $0$
can be arbitrarily slow.
\end{enumerate}
\end{theorem}

\noindent
{\bf Proof}: (a) Let $\alpha\in\R$ and $0<\beta-\alpha\leq 2\pi$.  Then,
by
Theorem~57 (p.~58) in \cite{celidze},
we can interchange the orders of integration to compute
$$
\int_{\alpha}^{\beta}u_r=
\int_{-\pi}^{\pi}f(\phi)v_r(\phi)\,d\phi,
$$
where 
$v_r(\theta)=P[\chi_{[\alpha,\beta]}](\ret)$. 

If $\beta-\alpha=2\pi$  then $v_r=1$  and the result is immediate.  Now
assume $0<\beta-\alpha<2\pi$.  For fixed $r$, the function $v_r$ has a maximum
at $\phi_1:=(\alpha+\beta)/2$
and a minimum at $\phi_2:= \phi_1+\pi$.  Use the Bonnet form of the
Second Mean Value Theorem for integrals (\cite{celidze}, p.~34) to write
\begin{eqnarray*}
\int_{\alpha}^{\beta}u_r
& = & 
\int_{\phi_1}^{\phi_2}f(\phi)v_r(\phi)\,d\phi
+\int_{\phi_2}^{\phi_1+2\pi}\!\!f(\phi)v_r(\phi)\,d\phi\\
& = &
v_r(\phi_1)\int_{\phi_1}^{\xi_1}f(\phi)\,d\phi+
v_r(\phi_1)\int_{\xi_2}^{\phi_1+2\pi}
\!\!f(\phi)\,d\phi\\
& = &
v_r(\phi_1)\int_{\xi_2-2\pi}^{\xi_1}\!\!f(\phi)\,d\phi
\end{eqnarray*}
where $\phi_1<\xi_1<\phi_2$ and $\phi_2<\xi_2<\phi_1+2\pi$.
And,
\begin{eqnarray*}
\left|\int_\alpha^\beta u_r\right| & \leq & 
\max_{\theta\in [-\pi,\pi]}v_r(\theta)
\left|\,\int_{\xi_2-2\pi}^{\xi_1}f\right|\\
& \leq & \|f\|.
\end{eqnarray*}
It follows that $\|u_r\|\leq \|f\|$.

(b) Let $\alpha\in\R$ and $0<\beta-\alpha\leq 2\pi$.
We have
\begin{eqnarray}
\int_{\alpha}^\beta \left[u_r(\theta)-f(\theta)\right]\,d\theta & = & 
\int_{\alpha}^\beta\left[\int_{-\pi}^\pi
\Phi_r(\phi-\theta)f(\phi)\,d\phi
-f(\theta)\int_{-\pi}^\pi
\Phi_r(\phi)\,d\phi\right]d\theta\notag\\
 & = & \int_{-\pi}^\pi\Phi_r(\phi)\int_{\alpha}^\beta
\left[f(\theta+\phi)-f(\theta)\right]\,d\theta\,d\phi.\label{5}
\end{eqnarray}
The reversal of integrals in \eqref{5} is justified by 
\cite[Theorem~58, p.~60]{celidze}.  We now have
\begin{eqnarray*}
\|u_r-f\| & \leq & \lsup_{0\leq\beta-\alpha\leq 2\pi}\left|\,
\int_{-\pi}^\pi\Phi_r(\phi)\int_{\alpha}^\beta
\left[f(\theta+\phi)-f(\theta)\right]\,d\theta\,d\phi\right|\\
 & \leq & P[g](r)\quad\text{where } g(\phi)=\|f(\phi+\cdot)-f(\cdot)\|.
\end{eqnarray*}
But if $f\in\hk$ then $f$ is continuous in the Alexiewicz norm, i.e.,
$\|f(\phi+\cdot)-f(\cdot)\|\to 0$ as $\phi\to 0$.  
See \cite{talvilacontinuity}.
Hence, $g$ is continuous at $0$ so $P[g](r)\to 0$ as $r\to 1$.

(c) Let $f$ be positive on $(0,1)$ and  vanish elsewhere.  Then $u_r$ is
positive for $0\leq r <1$.  We then have
\begin{eqnarray*}
\|u_r-f\| & \geq & \int_{-\pi}^0u_r(\phi)\,d\phi\\
 & = & \int_{0}^{1}f(\theta)P[\chi_{[-\pi,0]}](\ret)\,d\theta.
 \end{eqnarray*}

Now, as $r\to 1$
$$
P[\chi_{[-\pi,0]}](\ret)\to\left\{
	\begin{array}{cl}
	0, & 0<\theta<\pi\\
	1/2, & \theta=-\pi, 0, \pi\\
	1, & -\pi<\theta<0.
\end{array}
\right.
$$
But, the convergence is not uniform.
Let a decay rate be given by $A\fn [0,1]\to(0,1/2)$, where $A(r)$ decreases to
$0$ as $r$ increases to $1$. It is easy to show, for example, using a cubic
spline,
that $A$ has a decreasing $C^1$ majorant
with limit $0$ as $r\to 1$. So, we can assume $A\in C^1([0,1))$.
By keeping $\theta$ close enough to $0$ we can keep $P[\chi_{[-\pi,0]}](\ret)$
bounded away from $0$ for all $r$.
To see this, write $\rho:=(1+r)/(1-r)$.  Then
\begin{eqnarray*}
\|u_r-f\| & \geq &
\int_{0}^{1-r}\!\!\!\!f(\theta)P[\chi_{[-\pi,0]}](\ret)\,d\theta\\
 & = & \frac{1}{\pi}\int_{0}^{1-r}\!\!\!\!f(\theta)\left\{
	 \frac{\pi}{2}-\arctan\!\!\left[\rho\tan\!\!\left(\frac{\theta}{2}\right)
	 \right]+\arctan\!\!\left[\frac{1}{\rho}\tan\!\!\left(\frac{\theta}{2}\right)
	          \right]\right\}d\theta\\
 & \geq & \int_{0}^{1-r}\!\!\!\!f(\theta)\left\{
	          \frac{1}{2}-\frac{1}{\pi}
		  \arctan\!\!\left[\rho\tan\!\!\left(\frac{\theta}{2}\right)
                   \right]\right\}d\theta\\
 & \geq & \int_{0}^{1-r}\!\!\!\!f(\theta)\left\{
	 \frac{1}{2}-\frac{\rho\,\theta}{2\pi\cos(\theta/2)}\right\}d\theta\\
 & \geq & \left(\frac{1}{2}-\frac{1}{\pi\cos(1/2)}\right)\int_{0}^{1-r}\!\!\!\!f(\theta)
\,d\theta.
\end{eqnarray*}
We can now let
$$
f(\theta):=\left\{\begin{array}{cl}
-\left(\frac{1}{2}-\frac{1}{\pi\cos(1/2)}\right)^{-1}A'(1-\theta), 
& 0<\theta<1\\
0, & \text{ otherwise}.
\end{array}
\right.
$$
And,
$$
\|u_r-f\|\geq -\int_{0}^{1-r}A'(1-\theta)\,d\theta=A(r).\qed
$$
\begin{remarks}\label{remarks2}
{\rm
\begin{enumerate}
\item We have equality in (a) when $f$ is of one sign.
\item Part (a) and dilation show that if $0\leq r\leq s <1$ then
$\|u_r\|=\|P[u_s]_{\tfrac{r}{s}}\|\leq \|u_s\|$ (cf. \cite[Corollary~6.6]
{axler}).
\item The triangle
inequality and (b) show that $\|u_r\|\to\|f\|$ as $r\to 1$.
\item In (c),
$\|u_r-f\|$ can decay to $0$ arbitrarily fast.  Take $f$ to be constant!
\item The same proof shows that we can choose  $f\in L^p$ to make 
$\|u_r-f\|_p$ tend to $0$ arbitrarily slowly. For $1\leq p<\infty$,
let 
$$
f(\theta):=\left\{\begin{array}{cl}
\left(\frac{1}{2}-\frac{1}{\pi\cos(1/2)}\right)^{-1}p^{1/p}\left[A(1-\theta)
\right]^{1-1/p}
\left[-A'(1-\theta)\right]^{1/p}, 
& 0<\theta<1\\
0, & \text{ otherwise}.
\end{array}
\right.
$$
and then $\|u_r-f\|_p \geq A(r)$.
\qed
\end{enumerate}
}
\end{remarks}

\begin{theorem}\label{isometry}
The mapping
$P\fn\hk\to h^\hk$, $f\mapsto P[f]$, is an isometry into, but not onto,
$h^\hk$.
\end{theorem}

\noindent
{\bf Proof:}
Let $f\in\hk$ and $u=P[f]$.  From Remarks~\ref{remarks2}.2 and 
\ref{remarks2}.3,
$$
\|u\|_{\hk}=\lsup_{0\leq r<1}\|u_r\|=\llim_{r\to 1}\|u_r\|=\|f\|.
$$
Hence, $P$ is an isometry.

However, $P$ is not onto $h^\hk$.  Let $F$ be continuous on $[-\pi,\pi]$
such that $F(-\pi)=0$, $F$ is $2\pi$-periodic and $F$ is not in $ACG^*$, i.e.,
$F$ is not an indefinite Henstock--Kurzweil integral.  See \cite{celidze}
for the definition of $ACG^*$.  The function 
\begin{equation}
v_r(\theta):=F(\pi)\Phi_r(\pi-\theta)-\int_{-\pi}^\pi
\Phi'_r(\phi-\theta)F(\phi) \,d\phi\label{7.2}
\end{equation}
is harmonic in $D$ (using dominated
convergence).  Let $\alpha\in\R$ and $0<\beta-\alpha\leq 2\pi$.  Then
\begin{eqnarray*}
\int_{\alpha}^\beta v_r(\theta)\,d\theta & = & F(\pi)\int_\alpha^\beta
\Phi_r(\pi-\theta)\,d\theta-
\int_{-\pi}^\pi F(\phi)\int_{\alpha}^\beta
\Phi'_r(\phi-\theta)\,d\theta \,d\phi\\
 & = & F(\pi)P[\chi_{[\alpha,\beta]}](-r)
+ P[F](re^{i\alpha})-P[F](re^{i\beta}).
\end{eqnarray*}
So, $\|v_r\|\leq 3\max|F|$ and $v\in h^\hk$.
If there was $f\in\hk$ such that $v=P[f]$ then write 
$G(\theta):=\int_{-\pi}^\theta f$.  Since $G\in ACG^*$, we have
\begin{equation}
v(\ret)=G(\pi)\Phi_r(\pi-\theta)-\int_{-\pi}^\pi
\Phi'_r(\phi-\theta)G(\phi) \,d\phi.\label{7.3}
\end{equation}
Comparing \eqref{7.2} and \eqref{7.3}, letting $r\to 0$ shows $G(\pi)=F(\pi)$.
Write $H:=F-G$.  Expand $\Phi'_r(\theta)=
(-1/\pi)\sum_{n=1}^\infty nr^n\sin(n\theta)$.
The series converges uniformly and absolutely on compact subsets of $D$.  Then
for all $\ret\in D$,
\begin{eqnarray*}
0 & = & \int_{-\pi}^\pi H(\phi)\lsum_{n=1}^\infty nr^n\sin[n(\phi-\theta)]\,
d\phi\\
 & = &  \lsum_{n=1}^\infty nr^n\int_{-\pi}^\pi H(\phi)\sin[n(\phi-\theta)]\,
d\phi.
\end{eqnarray*}
For all $n\geq 1$ and all $\theta\in\R$ we have 
$\int_{-\pi}^\pi H(\phi)\sin[n(\phi-\theta)]\,
d\phi=0$.  Since $H$ is continuous it is constant.  But then $F$ differs
from $G$ by a constant.  This contradicts the assumption that $F\not\in ACG^*$.
Thus, no such $F$ exists and $P$ is not onto $h^\hk$.\qed

\section{Bounded variation}
Define the $2\pi$-periodic functions of normalised bounded
variation by
${\cal NBV}:=\{g\fn\R\to\R\mid g \text{ is } 2\pi\text{-periodic}, Vg<\infty,
g(-\pi)=0, g \text{ is right continuous}\}$.
Using the variation as a norm, ${\cal NBV}$ is a Banach space that is the
dual of $\hk$ \cite{swartz2}.  Analogues of Theorems~\ref{theorem1} and 
\ref{theorem2} now 
take the following form.
\begin{theorem}\label{bv}
Let $g\in{\cal BV}$
and
$v=P[g]$.
\begin{enumerate}
\item[{\rm (a)}] If $g\in {\cal NBV}$ then $v_r\to g$ weak* in ${\cal NBV}$ 
as $r\to 1$.
\item[{\rm (b)}] For all $0\leq r<1$, $\|v_r\|_\infty\leq
\inf|g|+Vg$.
\item[{\rm (c)}] If $g\in{\cal NBV}$ then
$\|v_r\|_\infty\leq
Vg$ for all $0\leq r<1$.
\item[{\rm (d)}] $Vv_r\leq Vg$ for all $0\leq r <1.$ 
\item[{\rm (e)}] There is $\sigma\in{\cal NBV}$ such that if
$w_r(\theta)=P[\sigma](\ret)$ then $V[w_r-\sigma]\not\to 0$
as $r\to 1$.  And, there is $\tau\in \bv$ such that if $w_r(\theta)
=P[\tau](\ret)$ and $\tau(\theta)=
[\tau(\theta+)+\tau(\theta-)]/2$ for all $\theta\in[-\pi,\pi]$ then
$V(w_r-\tau)\not\to 0$ as $r\to 1$.
\item[{\rm (f)}] Let 
$h^{\bv}:=\{u\fn D\to\R\mid\Delta u=0, \|u\|_\bv<\infty\}$,
where $\|u\|_{\bv}:=\sup_{0\leq r<1}Vu_r$.  The mapping
$P\fn{\cal NBV}\to h^\bv$, $g\mapsto P[g]$, is an isometric isomorphism
between the Banach spaces ${\cal NBV}$ and $h^\bv$.
\end{enumerate}
\end{theorem}

\noindent
{\bf Proof:} (a) Let $f\in \hk$.  Write $u=P[f]$.  Then, using \eqref{4.5} and
(b) of Theorem~\ref{theorem2},
\begin{eqnarray}
\left|\int_{-\pi}^\pi f(v_r-g)\right| & = &  
\left|\int_{-\pi}^\pi (u_r -f)g\right|\label{7.5}\\
 & \leq & \|u_r -f\| Vg\\
 & \to 0 & \text{ as } r\to 1.
\end{eqnarray}
The interchange of orders of integration in \eqref{7.5} is valid by
\cite[p.~58, Theorem~57]{celidze}.

(b), (c) These follow immediately from \eqref{4.5}. 

(d) Let $\{(s_n,t_n)\}$ be a sequence of disjoint intervals in $(-\pi,\pi)$.
Then
\begin{eqnarray*}
\sum|v_r(s_n)-v_r(t_n)| & = & \sum\left|\int_{-\pi}^\pi \Phi_r(\phi)
\left[g(\phi+s_n)-g(\phi+t_n)\right]d\phi\right|\\
 & \leq & P[1](r)\,Vg\\
 & = & Vg.
\end{eqnarray*}

(e) Let $-\pi<a<b<\pi$, $\sigma=\chi_{[a,b)}$ and 
$w_r(\theta)=P[\sigma](\ret)$.  Then $\sigma\in{\cal NBV}$ and
\begin{eqnarray*}
|w_r(b)-\sigma(b)-w_r(-\pi)+\sigma(-\pi)| & = & w_r(b)-w_r(-\pi) \\
 & \to & 1/2\quad\text{ as } r\to 1.
\end{eqnarray*}
So, $V(w_r-\sigma)\not\to 0$.

Note that if we replace $\tau(\theta)$ by
$[\sigma(\theta+)+\sigma(\theta-)]/2$ and now let
$w_r(\theta)=P[\tau](\ret)$ then $w_r(\theta)\to \tau(\theta)$
for all $\theta\in[-\pi,\pi]$.  But, 
$V(w_r-\tau)\to 2$ as $r\to 1$.  (Since $w_r(a)$ and $w_r(b)\to1/2$
as $r\to 1$.)

(f) Let $\sigma\in{\cal NBV}$ and $w_r(\theta)=P[\sigma](\ret)$.
By (d), $\|w\|_{\bv}\leq V\sigma$.  From (a), $w_r\to \sigma$ weak* in ${\cal NBV}$,
hence (cf. \cite[6.8]{axler}),
$$
V\sigma\leq \liminf_{r\to 1}Vw_r\leq \liminf_{r\to 1}\|w\|_{\bv}
=\|w\|_{\bv}.
$$
And, $P$ is an isometry.

To show $P$ is onto $h^{\bv}$, let $w\in h^\bv$.  Since $\hk$ is separable
\cite{swartz2},
every norm-bounded sequence in $\hk^*$ contains a weak* convergent subsequence
\cite[Theorem~6.12]{axler}.
But $\{w_r\}$ is norm-bounded in ${\cal NBV}$ so there is a subsequence
$\{w_{r_j}\}$  and $\sigma\in{\cal NBV}$ such that for all $f\in\hk$ we have
$\int_{-\pi}^\pi f w_{r_j}\to \int_{-\pi}^\pi f \sigma$ as $r_j\to 1$.  To show
$w=P[\sigma]$, fix $\ret\in D$. Then, since each function $w_{r_j}$ is
continuous on ${\overline D}$ and harmonic in $D$ it is the Poisson
integral of its boundary values, i.e.,
\begin{equation}
w(r_j\ret)=\int_{-\pi}^\pi \Phi_r(\phi-\theta)w_{r_j}(\phi)\,d\phi.\label{7.7}
\end{equation}
Now, $w$ is continuous on $D$, $\Phi_r(\cdot-\theta)\in\hk$ and $w_{r_j}$
is of bounded variation on $\partial D$, uniformly for $j\geq 1$.
Using weak* convergence, taking the limit
$r_j\to 1$ in \eqref{7.7} yields $w(\ret)=P[\sigma](\ret)$.  Thus, ${\cal NBV}$
and $h^\bv$ are isomorphic.  Since ${\cal NBV}$ is a Banach space, $h^\bv$ is
as well.\qed

\section{The Dirichlet problem}
Under an Alexiewicz norm boundary condition, we can prove uniqueness
for the Dirichlet problem.

\begin{theorem}\label{theorem4}
Let $f\in\hk$.  The Dirichlet problem 
\begin{align}
&u\in C^2(D)\label{8}\\
&\Delta u =0\quad \text{in } D\label{9}\\
&\|u_r-f\|\to 0 \text{ as } r\to 1\label{10}
\end{align}
has the unique solution $u=P[f]$.
\end{theorem}

\bigskip
\noindent
{\bf Proof:} First note that from Theorem~\ref{theorem2}(b) and
\cite[Proposition~1]{benedicks}, 
$u=P[f]$ is certainly a solution of \eqref{8},
\eqref{9} and \eqref{10}.

Suppose there were two solutions $u$ and $v$.  Write $w=u-v$.  Then
$w$ satisfies \eqref{8} and
\eqref{9}.  And, $\|w_r\|\leq \|u_r-f\|+\|v_r-f\|$, which has limit
$0$ as $r\to 1$.  Since $w$ is harmonic in $D$ it has the trigonometric
expansion
\begin{equation}
w(\ret)=\frac{a_0}{2}+\lsum_{n=1}^\infty r^n\left[a_n\cos(n\theta)+
b_n\sin(n\theta)\right],\label{11}
\end{equation}
the series converging uniformly and absolutely on compact subsets of $D$.
Fix $0\leq r<1$.  We have
$\|w_r\|\geq |\int_{-\pi}^\pi w_r|=\pi|a_0|$.  Letting $r\to 1$ shows
$a_0=0$.  And, for $n\geq 1$, we have $\|w_r\cos(n\,\cdot)\|\geq 
|\int_{-\pi}^\pi w_r(\theta)\cos(n\theta)\,d\theta|=\pi r^n|a_n|$.
As well, 
\begin{eqnarray*}
\|w_r\cos(n\,\cdot)\| & \leq & \|w_r\|\left\{\linf_{|\theta|\leq\pi}|\cos
(n\theta)|+V[\theta\mapsto\cos(n\theta)]\right\}\\
 & = & 4n\|w_r\|.
\end{eqnarray*}
Therefore, $4n\|w_r\|\geq \pi r^n|a_n|$.  Letting $r\to 1$ shows
$a_n=0$.  Similarly, $b_n=0$.  It follows that $w=0$ and we have
uniqueness.\qed

In \cite{shapiro},  Shapiro gave a uniqueness theorem that combined
a pointwise limit  with an $L^p$ condition.  There is an analogue for
the Alexiewicz norm.
\begin{theorem}\label{theorem5}
Suppose $\Delta u=0$ in $D$ and there exists $f\in\hk$ such that
\begin{align}
&u_r(\theta)\to f(\theta)\quad\text{ for each } \theta\in[-\pi,\pi)\label{12}\\
&\|u_r\|=o(1/(1-r))\quad\text{ as } r\to 1. \label{13}
\end{align}
Then $u=P[f]$.
\end{theorem}

\noindent
{\bf Proof:} As in Theorem~\ref{theorem4}, suppose $w$ is a solution of
the corresponding homogeneous problem ($f=0$).  Let $\alpha,\beta\in\R$
with $0<\beta-\alpha\leq 2\pi$.  Following the proof of Theorem~3 in
\cite{shapiro} and using \eqref{4.5},
\begin{eqnarray*}
\left|w(r^2 e^{i\theta})\right|
 & = & \left|P[w_r](\ret)\right|\\
 & \leq & \frac{\|w_r\|\,g(r)}{2\pi(1-r)},
\end{eqnarray*}
where $g(r):=(1+6r+r^2)/(1+r)$.  But, $g(r)\leq g(1)=4$.  Hence, by 
\eqref{13},
$\|w_{r^2}\|_\infty=o(1/(1-r)^2)$ and so 
$\|w_r\|_\infty=o(1/(1-r)^2)$ as $r\to 1$.  
It follows
from \cite[Theorem~1]{shapiro} that $w=0$.\qed

As pointed out in \cite{shapiro}, neither \eqref{12} nor \eqref{13} can be
relaxed.  If $u_r\to f$ except for one value $\theta_0\in[-\pi,\pi)$
then we can add a multiple of $\Phi_r(\theta -\theta_0)$ to $u(\ret)$.
If in place of \eqref{13} we have $\|u_r\|=O(1/(1-r))$ then we can
add a multiple of $\Phi_r'$ to $u_r$,
since for each $\theta\in\R$, $\Phi'_r(\theta)\to 0$ as $r\to 1$.

\begin{example}\label{example}
{\rm
(a) Let $f\in\hk\setminus L^1$.  Then the unique solution to 
\eqref{8}--\eqref{10} is $u=P[f]$.  In this case, the $L^p$ norms
of
$u_r$ need not be bounded as $r\to 1$.  If we are given a harmonic
function $v$ such that the Alexiewicz norms $\|v_r\|$ are uniformly 
bounded for $0\leq r<1$ 
then we cannot infer the existence of $g\in\hk$ such that
$v=P[g]$.  This is because $\hk$ is not complete.

(b) Let $v(z)=(1+z)/(1-z)$ and $w(z)=v(z)e^{-v(z)}$. Define
\begin{eqnarray*}
u(\ret)  & = & Re(w(\ret))  \\
 & = & \frac{(1-r^2)\cos\left(
\frac{
2r\sin\theta}{1-2r\cos\theta+r^2}\right)+2r\sin\theta\sin\left(
\frac{
2r\sin\theta}{1-2r\cos\theta+r^2}\right)}{
\exp(2\pi\Phi_r(\theta))
(1-2r\cos\theta+r^2)}.
\end{eqnarray*}
Let 
$$
f(\theta):=\llim_{r\to 1}u_r(\theta)=\left\{
\begin{array}{cl}
\left(\frac{\sin\theta}{1-\cos\theta}\right)
\sin\left(\frac{\sin\theta}{1-\cos\theta}
\right), & 0<|\theta|<\pi\\
0, & |\theta|=0, \pi.
\end{array}
\right.
$$
Note that $f\not\in L^p$ for any $1\leq p\leq \infty$.
The set function $\mu$ defined by $\mu(A)=\int_A f$ is not a signed Borel
measure.  Thus, $u$ is not the Lebesgue--Poisson integral of any
$L^p$ function or measure.  Since $f(\theta)\sim (2/\theta)\sin(2/\theta)$
as $\theta\to 0$ we have $f\in\hk$.  And,
\begin{eqnarray*}
|(1-r)u_r(\theta)| & \leq & (1-r)e^{-1}+\frac{2r\,e^{-1}}{1+r}\\
 & \leq & 1/2.
\end{eqnarray*}
By dominated convergence, $\|(1-r)u_r\|\to 0$ as $r\to 1$.  And,
by Theorem~\ref{theorem5}, $u=P[f]$.  There is a similar result for
the imaginary part of $w$.

(c) Let $w(z)=[1/(1-z)]e^{[1/(1-z)]}$ and define
\begin{eqnarray*}
u(\ret)  & = & Re(w(\ret))  \\
 & = & \frac{
(1-r\cos\theta)\cos\left(
\frac{
r\sin\theta}{1-2r\cos\theta+r^2}\right)-r\sin\theta\sin\left(
\frac{
r\sin\theta}{1-2r\cos\theta+r^2}\right)}{
\exp\left(\frac{r\cos\theta-1}{1-2r\cos\theta+r^2}\right)
(1-2r\cos\theta+r^2)}.
\end{eqnarray*}
Let 
$$
f(\theta):=\llim_{r\to 1}u_r(\theta)=\left\{
\begin{array}{cl}
\sqrt{e}\left[
\frac{(1-\cos\theta)\cos\left(\frac{\sin\theta}
{2(1-
\cos\theta)}\right) -
\sin\theta\sin\left(\frac{\sin\theta}
{2(1-\cos\theta)}\right)}{2(1-\cos\theta)}\right]
, & 0<|\theta|\leq\pi\\
\infty, & \theta=0.
\end{array}
\right.
$$
Although $f\in \hk$,
Theorem~\ref{theorem5} does not apply since $f$ is not a real-valued
function.  Indeed, $(1-r)u_r(0)=\exp(1/(1-r))\to\infty$ as $r\to 1$.
From Theorem~\ref{theorem1}, $u$ is not the Poisson integral of any
function in $\hk$ (nor $L^p$ function nor measure).  In particular,
$u\not=P[f]$.
}
\end{example}

In examples (b) and (c), the origin is the only point of nonabsolute summability
of $f$.  For each $0\leq\lambda<2\pi$, an example  is given in \cite{benedicks}
of the Poisson integral
of a function in $\hk$
whose set of points of nonabsolute summability in $(-\pi,
\pi)$ has measure $\lambda$.

\end{document}